# MOUVEMENT
# &
# ORIGINE DU CALCUL INFINITÉSIMAL :

### Théorisation du mouvement et infinitésimaux

Salomon OFMAN
Institut mathématique de Jussieu-Paris Rive Gauche
Histoire des Sciences mathématiques
4 Place Jussieu
75005 Paris
salomon.ofman@imj-prg.fr

*English Abstract.* – In a previous article ([OFM2]) we gave the general foundations of the theory of movement considered from a philosophical and mathematical point of view. Philosophical it meant to understand the opposition of the one and the multiple, mathematically to consider the opposition between the discreet and the continuous.
In this article we want to show how the widespread introduction of mathematics in physics by Galilei leas to a change in the very notion of movement. Conversely, this theorization of the movement is the origin of the elaboration of the infinitesimals, one of the major upheavals in mathematics.

Résumé. – Dans un article précédent ([OFM2]), nous avons posé les fondements généraux de la théorie du mouvement d'un point de vue philosophique et mathématique. Sur le plan philosophique, il s'agissait de comprendre l'opposition de l'un et du multiple qui, sur le plan mathématique, prenait l'aspect de l'opposition du discret et du continu.
Dans ce présent travail, nous nous proposons de montrer comment l'introduction massive par Galilée des mathématiques dans la physique, conduit à un changement de la notion même de mouvement. Inversement, cette théorisation physique du mouvement est à la source de l'élaboration des infinitésimaux, l'un des bouleversements majeurs en mathématiques.

Toutes les références aux œuvres de Galilée renvoient à l'Édition Nationale italienne.

## I. Aristote et la physique du mouvement.

### 1. Certains concepts de la *Physique* d'Aristote.

Nous allons exposer brièvement la théorie du philosophe sur le mouvement. Une difficulté, récurrente pour l'analyse de la pensée aristotélicienne, est qu'on ne trouve pas toujours un développement univoque dans ses différents ouvrages scientifique qui, d'une certaine manière, traitent tous de cette question. Pour Aristote, en effet, toute donnée physique est sujette au mouvement.

Notre démarche adopte un sens téléologique, l'objectif principal étant de souligner les traits qui se retrouvent au fondement de la conception du mouvement jusqu'à l'époque médiévale. C'est explicitement contre cette théorie (cf. la lettre de Galilée à Fortunio Liceti, 15 septembre 1640) que sera fondée la physique moderne, qui débute avec Galilée.

Cette démarche téléologique est conforme, si ce n'est à l'histoire des idées, du moins à la philosophie du Stagirite. Sa conception de la causalité admet en effet quatre formes de causes, dont la principale, la cause finale, est ce en vue de quoi quelque chose se produit.

Les cinq points fondamentaux, quoique non exhaustifs, sur lesquels s'appuie la théorie du mouvement dans la *Physique* sont les suivants.

> A. *Le continu est ce qui est divisible sans fin (VI, 1, 231b15).*
> B. *Aucun infini n'existe réellement ou en réalité[1].*
> C. *La Physique (i.e. la Nature) est un principe de mouvement et de repos (κινεῖσθαι καὶ ἠρεμεῖν) (II, 1, 192b21-22).*
> D. *Physiciens et mathématiciens diffèrent en ce que les premiers étudient (enquêtent sur) les corps naturels, les seconds (sur) des abstractions (points, lignes, plans, formes solides) issues de ces corps.*
> E. *Les catégories (i.e. les propriétés essentielles des êtres) de quantité et de qualité ne peuvent se changer l'une en l'autre. On comprend cela généralement comme l'affirmation d'une disjonction complète entre les deux. Sa conséquence la plus importante est négative, l'impossibilité pour ce qui appartient à l'une de ces catégories d'appartenir à une autre. Il faut toutefois noter que les termes mêmes 'quantité' et 'qualité', en tant que catégories, ne se trouvent pas dans les textes d'Aristote.*

Pour éviter les problèmes posés par les divers sens contemporains de 'vitesse', on utilisera dans la suite le terme 'célérité' pour traduire les multiples notions de rapidité d'un mouvement.

Il ne faudrait pas croire que ces ambiguïtés dérivent exclusivement de la langue. Par exemple, le terme classique latin de *velocitas* utilisé en scolastique pose les mêmes difficultés. Or suivant la traduction choisie, on aboutit à des interprétations très différentes.

Encore de nos jours, le concept de 'vitesse' est réduit aux cas les plus simples. En général, les dictionnaires la définissent simplement comme caractérisant la rapidité. Pour le Littré, c'est aussi la distance effectuée par unité de temps (Littré, 1958). Le *Trésor de la Langue Française du CNRS*, ajoute un sens technique, en mathématique et en physique, *'l'espace parcouru par un mobile dans l'unité de temps choisie ; [le] rapport entre ces deux éléments'*.

Il y a là une difficulté inhérente à ce type de notions, où il s'agit d'établir une relation entre une conception immédiate et sa définition physico-mathématique. Ainsi que le remarque Poincaré à propos du calcul infinitésimal, les questions qui se posent à propos d'un énoncé mettent généralement en évidence des difficultés liées à la nature de la chose dont on parle ([POI], p. 129). C'est le problème où, selon Platon, l'intelligence est en lutte avec elle-même (cf. [OFM3], §III.3).

---

[1] Le statut de l'infini aristotélicien est extrêmement complexe et nous ne pouvons l'aborder ici. Pour simplifier à l'extrême, situé entre l'être et le non-être, son *degré* de réalité est particulièrement bas. Cette question est étudiée plus en détail dans [OFM3], au travers d'une analyse du continu aristotélicien et de ses rapports à la topologie moderne.

## 2. La notion de mouvement chez Aristote.

La *Physique* d'Aristote ('περὶ φύσεως', qu'on devrait traduire plus précisément par 'enquête sur la nature') s'occupe des objets du changement (III, 1, 200b12). Le terme grec général pour 'changement' est μεταβολή, et cependant, celui utilisé la plupart du temps par Aristote est κίνησις (mouvement). Il faut donc comprendre que le mouvement, selon Aristote, concerne tout changement et pas seulement celui d'un corps mobile (ainsi le chaud et le froid, la moralité et l'amoralité, la justice et l'injustice). La vélocité célérité du mouvement étant une qualité de ce mouvement (en français/italien, on parle aussi d'un mouvement, d'une voiture…, rapide ou véloce), on conclut de la propriété 1.E qu'elle ne peut être mesurée.

Pourtant un mobile peut être plus ou moins rapide qu'un autre. Une qualité, selon Aristote, est en effet sujette au plus et au moins, un homme plus ou moins heureux, un manteau plus ou moins blanc et naturellement Achille plus rapide que la (tortue) plus lente. La disjonction aristotélicienne qualité/quantité se comprend donc non du point de vue de l'addition/soustraction mais de la multiplication/division. Il est possible d'accroître ou de diminuer une qualité comme une quantité. À l'inverse, si l'on peut multiplier et diviser la seconde, il en va différemment de la première. Ainsi la possibilité de parler de plus rapide ou de plus lent ne pose guère de problème du point de vue d'une qualité. Une quantité est mesurable parce qu'il est une certaine unité (IV, 12, 220b14-24). Mais aucune unité n'existe pour une qualité, ainsi la justice et l'injustice, le blanc et le non blanc, le chaud et le froid et aussi bien la célérité et la lenteur.

Associé à 'l'interdit canonique' de division des quantités hétérogènes pesant sur les sciences grecques, le système aristotélicien aurait ainsi rendu impossible toute mathématisation du mouvement. Celle-ci requérant la division d'une distance par un temps. C'est ce que nous appellerons par la suite l'interprétation standard.

## 3. Les problèmes que pose cette interprétation.

Pourtant, on trouve dans la littérature exégétique, diverses descriptions mathématiques du mouvement et même de la célérité. Ainsi, dans la classique traduction (et commentaires) de Ross (p. 28), celui-ci a recours, pour rendre compte du texte d'Aristote, aux formules suivantes :

**D = C•M•T/Densité du milieu      (1)**

**V = C•M/ Densité du milieu      (2)**

où le milieu est ce en quoi le mobile se meut (air, eau…), M la masse du mobile, T le temps et C une constante (le signe '•' entre C, M et T est un ajout, dans le commentaire original le signe multiplicatif n'est pas indiqué).

Il faut noter que l'on a établi un rapport entre des termes supposés être purement qualitatifs. Or ces équations contredisent l'interdit numérique qui, selon l'interprétation standard, devrait s'attacher au mouvement. En effet, une première conséquence de ces formules, est la mesurabilité de la célérité, une qualité, en contradiction avec l'interprétation usuelle du point fondamental 1.E ci-dessus.

En outre, l'équation (1) n'aurait pas de sens, puisqu'elle utilise un rapport entre choses hétérogènes (masse, temps, densité). Selon Aristote, le temps est le nombre du mouvement (IV, 11, 219b2). Cette définition est très problématique et a été objet de critique depuis l'Antiquité grecque. Ainsi d'après Simplicius, déjà Straton de Lampsaque (autour de -300) objectait que le nombre est une quantité discrète (διωρισμένον), alors que le temps étant continu n'est pas nombrable (ἀριτμητόν). Mais du moins, tous les termes de (1) sont-ils des quantités.

Au contraire, l'équation (2) présuppose, elle, la possibilité d'une mesure mathématique de la célérité/rapidité, qui n'est pas une quantité, mais une qualité. Or une telle mesure n'apparaît pas dans le texte d'Aristote, et, toujours suivant l'interprétation standard, elle contredirait également 1.E.

Naturellement, et Ross l'écrit clairement, Aristote n'a jamais donné ces équations. Toutefois, le Stagirite affirme bien que le rapport (λόγος) des distances parcourues par deux mobiles de même célérité dans le même temps, est égal à l'inverse de la densité du milieu dans lequel ils se meuvent (IV, 8, 215a29-b12). Il est vrai qu'il ne considère que des rapports entiers. Mais encore Galilée utilise des rapports entiers dans ses exemples. On peut donc penser que c'est par simple commodité que l'un et l'autre s'y restreignent. Ainsi ces formules paraissent bien refléter la pensée d'Aristote concernant la célérité, et être le moyen le plus efficace pour la comprendre.

Et on trouve précisément dans la *Physique* une mise en rapport explicite de célérités, qui sont donc traitées, au sens euclidien, comme des grandeurs (cf. §I.3') :

'*Il y a donc la même proportion (λόγον) entre [la facilité de pénétration dans] l'air et [celle dans] l'eau qu'entre la vitesse (τάχος) dans l'un et la vitesse dans l'autre*' ([ARI], IV, 215b5-10),

ce qu'on peut réécrire en symboles mathématiques :

$C_1/C_2 = R_2/R_1$

où, pour un mobile donné, $C_1$ et $C_2$ sont respectivement sa célérité dans l'air et dans l'eau, $R_1$ et $R_2$ représentant quelque chose comme les résistances à ce mobile de l'air et de l'eau.

Il paraît donc très exagéré de soutenir que cette physique est radicalement incompatible avec toute interprétation mathématique (cf. [KOY3], p. 17). Ce qui ne signifie évidemment pas, qu'on ait, comme chez Galilée, une théorie mathématique vérifiable pas-à-pas par expériences.

À ces textes tirés du corpus aristotélicien, nous voudrions en ajouter deux autres plus anciens, rarement cités dans ce contexte. Dans *Le Politique*, Platon divise 'la technique de la mesure' en quatre parties : 'le nombre, la longueur, la profondeur, la largeur et la vitesse (ταχυτής) ([PLA3], 284e). Dans *Lachès*, Socrate prend comme modèle d'une définition, celle de la 'vitesse' (ταχυτής), et la définit en termes quantitatifs : 'Si donc quelqu'un me demandait : « Dis, Socrate, en quoi consiste ce que tu appelle 'vitesse' dans tous les cas ? », je lui répondrais qu'en ce qui me concerne, j'appelle vitesse la capacité de faire **plusieurs** chose en **peu** de temps' ([PLA4], 192a ; nous soulignons).

Il ne paraît donc nullement illégitime à Aristote, et plus généralement aux Grecs anciens, si ce n'est de multiplier et de diviser des termes hétérogènes, distances par temps, du moins de multiplier et diviser numériquement la célérité.

L'interprétation standard ne semble donc pas pouvoir rendre compte de manière satisfaisante des approches quantitatives de la *Physique*. Dans ce qui suit, nous allons voir qu'elle conduit à des difficultés plus sérieuses encore, à commencer dans les *Éléments* d'Euclide.

### 3'. Catégories aristotéliciennes et les *Éléments* d'Euclide.

Nous voudrions montrer dans ce paragraphe que la disjonction entre les diverses catégories, et en particulier entre qualité et quantité, n'est pas toujours respectée dans les textes mathématiques. Et pas même cette autre frontière, bien plus importante pour les mathématiciens, qui sépare le discret du continu ou l'arithmétique de la géométrie.

C'est au livre I, qu'Euclide traite des angles. Les mathématiciens grecs anciens considéraient non seulement les angles entre des droites, ce que nous appelons encore angles, mais entre des courbes, ainsi entre un cercle et une tangente à ce cercle. L'angle y est défini comme une qualité ou une relation (une inclinaison). Et pourtant, il y est traité comme une grandeur, en particulier il est toujours divisible. D'où l'interprétation de Proclus (un millénaire après) pour qui l'angle participerait des trois catégories.

Le livre V, quant à lui, construit une théorie extrêmement élaborée des proportions pour les grandeurs ou quantités (μεγέθει).

On y trouve une définition très générale des rapports et des égalités de rapport, valables pour toute grandeur. Toutefois, la grandeur elle-même n'est pas objet de définition ; elle admet un statut qu'on pourrait rapprocher (en mathématique contemporaine) de celui d'ensemble.

La théorie moderne des ensembles présente en effet la même difficulté définitionnelle. Ainsi, l'inclusion entre ensembles donne un ordre (non total). Cet ordre devrait toutefois être défini sur la totalité des ensembles. Ce qui pose le problème de cette totalité et conduit nécessairement à sortir des ensembles, afin d'éviter le paradoxe de Russell de 'l'ensemble' de tous les ensembles.

Rester dans un certain flou concernant la grandeur, permettait à la fois de maintenir la rigueur de la théorie et son caractère opératoire, tout en conservant sa très grande généralité. Ne pas définir (au sens habituel) la grandeur, loin d'être un manque, était donc une nécessité pour le mathématicien grec. C'est de manière très proche, ce qu'on retrouve dans la théorie mathématique moderne. Ainsi, dans son petit livre classique sur les bases de la théorie ensembliste, Paul Halmos affirme-t-il :

'*Il est une chose que le développement ne contiendra pas, c'est une définition des ensembles. La situation est analogue à l'approche axiomatique familière de la géométrie élémentaire. Cette approche **n'offre pas de définition des points et des lignes ; en place, elle décrit ce que l'on peut faire avec ces objets**'* (nous soulignons) ([HAL], p. 2).

La définition de l'ensemble y est remplacée par des règles d'utilisation. De même doit-on comprendre la notion de grandeur d'après son aspect opératoire. Ce sera essentiellement ce qui peut être d'une part **multipliée par un entier**, et d'autre part **comparée** (à d'autres grandeurs de même nature ou genre (μεγετῶν ὁμογενῶν)), suivant les. Cela suit des définitions V.4 et V.5, d'où leur puissance opératoire. Cette dernière définition est d'autant plus importante pour nous, que Galilée l'utilise pour établir sa propre définition de l'égalité des 'vitesses', plus précisément pour les quantifier (cf. §III.1).

Définition V.4 : '*Des grandeurs sont dites avoir un rapport l'une relativement à l'autre quand elles sont capables, étant multipliées, de se dépasser l'une l'autre*' ([EUC], p. 39).

Définition V.5 : '*Des grandeurs sont dites être dans le même rapport, une première relativement à une deuxième et une troisième relativement à une quatrième quand des équimultiples de la première et de la troisième ou simultanément dépassent, ou sont* *simultanément égaux ou simultanément inférieurs à des équimultiples de la deuxième et de la quatrième, selon n'importe quelle multiplication, chacun à chacun, [et] pris de manière correspondante.*' (ib.).

Ce dernier énoncé peut être réécrit, en termes modernes, de la manière suivante :
Pour tous entiers (positifs, distincts, non nuls) *m* et *n*, on a l'équivalence suivante :
[*E/F* = *G/H*] équivaut à [((*E/F* > *n/m*) ⇔ (*G/H* > *n/m*)) et ((*E/F* = *n/m*) ⇔ (*G/H* = *n/m*)) et ((*E/F* < *n/m*) ⇔ (*G/H* < *n/m*))],
la dernière condition pouvant être supprimée, si à l'instar de la plupart des commentateurs, on suppose qu'il existe un ordre total sur les grandeurs.

Ceux-ci discutent pour savoir si, dans les *Éléments*, les entiers peuvent être considérés comme des grandeurs, et si les rapports de grandeurs peuvent être traités comme des grandeurs. Cependant, ils s'accordent pour voir dans ces quelques définitions, les prémisses de la théorie moderne des réels. Dans ce cadre, l'équivalence ci-dessus se traduit en disant que l'ensemble des rationnels est (isomorphe à une partie) dense dans l'ensemble des rapports de grandeurs (une grandeur n'est pas nécessairement un nombre). En fait, la théorie euclidienne des grandeurs va bien au-delà, ainsi que le montre la possibilité de traiter les angles comme des grandeurs (cf. [OFM3], Annexe 6).

Quoiqu'il en soit, au livre X, les nombres entiers sont traités comme des grandeurs, ce qui a provoqué d'importantes polémiques sur le caractère rigoureux de ses propositions (prop. X.1 à X.11). Car Euclide transgresse ici, ce qui est considéré comme l'un des 'interdits' de la pensée grecque ancienne, la séparation entre le discret et le continu. Tant il est vrai, comme le note Leibniz, que pour un mathématicien, il n'est pas question d'entraver 'l'art d'inventer par excès de scrupules' ou 'rejeter sous ce prétexte les meilleurs découvertes, en nous privant de leurs avantages' (Leibniz, *Réponse à Nieuwentift*, [LEI1], V, p. 320-328).
Mais dès lors qu'une grandeur n'est rien d'autre que ce qui, d'une part peut être multipliée (par des entiers), et d'autre part être comparable à d'autres grandeurs de même nature, les entiers tout comme les célérités, qui suivant les textes aristotéliciens vérifie ces propriétés, peuvent être considérées comme des grandeurs.
Et comme on l'a vu, les angles sont traités dans les *Éléments* à la fois ou alternativement comme des qualités, des quantités ou des relations, transgressant une fois encore la division des catégories aristotéliciennes.

Il est donc difficile d'affirmer que les mathématiciens grecs se soient interdit l'étude du mouvement par conformité aux écrits du Stagirite, qui d'ailleurs n'avaient sans doute pas l'importance qu'ils acquièrent à l'époque médiévale (cf. par exemple [CR-PE], p. 238). De plus, il n'est pas impossible de concilier d'une certaine manière les interdits aristotéliciens et une quantification du mouvement. Les catégories en effet ne visent pas tant les mathématiques, dont les objets ne sont, pour lui, qu'une abstraction des choses, que les sciences s'occupant des choses elles-mêmes.
Au livre VI de la *Physique*, il utilise la célérité pour produire une relation entre le temps et la distance, prouvant par exemple la continuité et même l'existence (en un certain sens) du temps (232b20-233a20). Pour cela, il ramène le problème au rapport entre temps et célérité d'un mouvement uniforme. En effet raisonne-t-il, parce que l'on peut trouver un mobile aussi lent que l'on veut, le temps doit être indéfiniment divisible. Il établit ainsi un parallélisme complet, un 'isomorphisme', entre temps et mouvement. Si le temps nombre le mouvement, le mouvement nombre le temps. En termes modernes, on parlerait d'un isomorphisme entre ensembles ou d'ensembles de même cardinal. Aristote donne comme analogue du rapport entre temps et mouvement, celui entre 'un nombre' et 'le nombre de chevaux' ([ARI], IV, 220b14-20). Comme nous l'avons vu au paragraphe précédent, cette numérisation du temps et du mouvement est incompatible avec une séparation absolue entre discret et continu. On retrouve la situation du livre X des *Éléments*.

Dans la mesure des textes qui nous sont parvenus, ces témoignages nous paraissent montrer que la conception quantitative du mouvement s'inscrit de manière cohérente à l'intérieur du corpus scientifique de la Grèce ancienne, de Platon à Euclide, et s'accorde avec une interprétation opératoire de la grandeur en mathématiques. En n'excluant pas de traiter la célérité comme une grandeur, Aristote restait ainsi dans le cadre des sciences grecques.

Une question se pose alors. Si le problème pour les mathématiciens, et Aristote lui-même, n'était pas l'appartenance de la célérité à la qualité, comment comprendre l'absence de théorie mathématique du mouvement en Grèce ancienne ?

La difficulté nous paraît se trouver non pas dans l'alternative quantité/qualité, mais dans la nécessité d'un étalon pour mesurer la vitesse. Tout d'abord, il n'est pas facile d'obtenir un mouvement toujours identique servant de norme, de même que pour les longueurs (la règle), le poids ou même le temps (horloge à eau, sablier, pouls). Mais en outre, et c'est le problème le plus insoluble, il est une infinité de mouvements tous différents, soit par la trajectoire, soit par leurs vitesses. L'obstacle à toute étude quantitative était bien plutôt la multiplicité des vitesses (dès lors que la célérité n'est pas uniforme) et leur caractère infinitésimal. Une définition formelle n'est en effet possible que, d'une manière ou d'une autre, sous forme infinitésimale. Les autres conceptions se ramènent à celle de la vitesse moyenne. Celle-ci réduit tous les mouvements au mouvement uniforme qui n'est rien d'autre qu'un couple formé d'une distance (celle parcourue par le mobile) et d'un temps (la durée du mouvement) (cf. §III.3). Dans ce cadre, l'intérêt d'une dynamique paraît très limité.

Toutefois, à l'encontre de nous autres modernes, le mouvement pour Aristote, n'est qu'une partie d'une conception philosophique, dont l'objectif est de rendre compte de la totalité de l'univers, au sens grec d'ordre universel (le 'cosmos')[2]. L'erreur serait de le concevoir selon un sens contemporain hérité de la physique classique.

Il n'empêche que, pour nous aussi, la célérité a un double sens (cf. §I.3). L'un est qualitatif, l'autre quantitatif, ce qui ajoute à l'ambiguïté de cette notion, et rend plus difficile encore de saisir son sens passé.

La théorie du mouvement, élaborée par Aristote, devait permettre d'éviter les paradoxes, tels ceux de Zénon. Ainsi, considérer la célérité comme qualité et non quantité, justifiait que le mouvement fût objet d'étude de la Physique, tout en rejetant les paradoxes sur l'infini ('Achille et la tortue, 'la flèche'…). Dès lors, on comprend que les difficultés inhérentes à cette étude, induites par le manque d'unité du mouvement, aient pu bloquer sa mathématisation, bien plutôt que l'inverse.

En obtenant cette unité, la physique galiléenne va permettre d'élaborer une théorie mathématique du mouvement. On peut donc bien dire que le problème du mouvement est celui de son hétérogénéité. Mais non pas mathématique (diviser des termes hétérogènes), mais de sens : le concept de célérité étant en soi hétérogène.

## II. La question du mouvement jusqu'à Galilée.

---

[2] C'est de ce point de vue que Geymonat conçoit la défaite de ce qu'il appelle 'le programme de politique culturelle de Galilée' ([GEY], p. 78). La condamnation de Galilée par l'Église est l'achèvement d'une lutte interne entre aristotéliciens et modernistes (*ib.*, p. 113-114).

Comme cela était prévisible (cf. [OFM2], §III.1), deux écoles s'affrontent à propos de la naissance de cette nouvelle physique.

## 1. L'interprétation 'catastrophique'.

Pour simplifier, il n'y aurait eu pour ainsi dire aucun véritable progrès de l'Antiquité tardive (datée autour de -200 i.e. Archimède) à la Renaissance européenne (cf. par exemple [WHI], p. 17, 22-23).

À strictement parler, cela est évidemment inexact, même si l'on s'en tient aux mathématiques, beaucoup de travaux ayant eu lieu au Moyen-âge. Ainsi ceux qui aboutissent au développement de l'algèbre, ou encore ceux des savants d'Oxford et de Paris, portant sur les suites infinies et les surfaces (cf. §II.2). Et c'est bien avant la naissance du savant italien que sont introduits les étranges 'imaginaires' (i.e. les nombres complexes) pour résoudre certaines équations polynomiales. Mais il est vrai qu'on ne voit guère de changement brutal, ce qu'on pourrait nommer une révolution intellectuelle. Et en ce qui concerne la géométrie, les textes de base restent ceux d'Euclide, parfois mal traduits et mal compris.

On pourra penser que chercher la solution d'une question dans un 'génie', est une solution de facilité. Toutefois, cette interprétation, comme toutes celles qui peuvent être qualifiées de catastrophiques, a un caractère relatif. Dans ce cadre interprétatif, la disparition des dinosaures s'est étalée sur plusieurs dizaines de milliers d'années (cf. [OFM3], Annexes 1 à 3). De même, Galilée élabore sa nouvelle science du mouvement dans la durée, en relation avec les travaux des autres mathématiciens et expérimentateurs. Il suffit pour cela de considérer l'évolution de ses idées sur le sujet, mais également l'importance de la correspondance qu'il entretint. Il en va de même pour l'intérêt qu'il montre envers les découvertes techniques.

On en trouve un exemple avec la lunette astronomique qui lui permettra des découvertes cosmologiques si importantes. C'est par elle qu'il atteindra un statut exceptionnel et un prestige qui lui donneront une grande liberté de manœuvre. Toutefois, l'utilisation de verres correcteurs est déjà ancienne, puisqu'on la trouve sur une illustration datant de la fin du 15$^{ème}$ siècle, la *Nef des Fous* de Sebastian Brant (1494) (NY Academy of Medicine). Quant à l'invention de la lunette astronomique, vers 1608, elle est attribuée à un opticien hollandais Hans Lippershey. Et c'est par d'autres scientifiques, que Galilée en aura connaissance, ce qui lui permit du reste d'en construire lui-même.

Les commentateurs divergent sur l'interprétation de l'œuvre de Galilée, souvent en raison de l'évaluation qu'ils en font, et de la place qu'ils lui attribuent dans la science. Certains comme A. Koyré privilégient son platonisme et sa méthode théorique, d'autres comme L. Geymonat insistent au contraire sur ses expériences voire son empirisme expérimental. Toutefois, ce qui importe ici, est l'étude du mouvement qui se développe à cette époque et va conduire à une élaboration de la mécanique et des infinitésimaux. Ce n'est donc pas l'originalité des travaux galiléens, en rapport à ses contemporains ou prédécesseurs proches, qui nous occupe, mais leur impact en tant qu'ils débouchent sur une nouvelle conception de la physique. Notre propos en effet, n'est pas d'ajouter un nom à la liste déjà longue des proto-découvreurs du calcul infinitésimal, mais de montrer en quoi le travail de Galilée a été à son fondement.

C'est en ce sens, et non au travers de la problématique du génie en histoire/science, que nous comprenons l'interprétation 'catastrophiste'.

## 2. Le point de vue 'évolutionniste' et la question de l'infini.

Pour les Scolastiques, la séparation ou non-séparation des catégories de la qualité et de la quantité posait un problème 'philosophique' fondamental. Par contre rien ne prouve qu'il en fût ainsi dans l'Antiquité grecque. Platon, analysant les grand systèmes politiques, dans ce que les interprètes modernes considèrent souvent comme une plaisanterie, n'hésite pas à quantifier le rapport de bonheur entre le roi/philosophe et le tyran, précisément égal à $(3^2)^3 = 3^6 = 9^3 = 729 \approx 2 \times 365$ (1 journée heureuse de l'un équivaut approximativement à 2 années de l'autre) (*République*, IX, 587d-588a).

Selon Thomas d'Aquin, le plus ou moins de qualité s'explique par un processus d'augmentation/diminution de participation à cette qualité ; pour Durant de Saint-Pourçain (13$^{ème}$ siècle), cela passe par l'accès à différents degrés de perfection (devenir plus chaud signifiant être plus parfait relativement à la chaleur).

Pour Duns Scot et son école des Scotistes (début du 14$^{ème}$ siècle), cette variation de qualité est causée par l'addition/soustraction de différents degrés (de qualité) qui s'additionnent comme des gouttelettes d'eau. Un mouvement peut alors prendre successivement des degrés de célérité (toujours une qualité) variables, et il est tentant de comprendre leur augmentation ou leur diminution comme une sorte d'accélération/décélération. Cette conception est adoptée par Guillaume d'Ockham ou Occam, puis par la plupart des 'philosophes' scolastiques.

A peu près à la même époque, les physiciens d'Oxford et de Paris étudiaient ces mêmes problèmes en relation avec l'infini mathématique, en particulier au travers de séries infinies, et les appliquaient au problème de la célérité.

Le mouvement d'un mobile était alors classifié par son 'uniformité' ou sa 'difformité'. On dit qu'un mouvement est uniforme (respectivement difforme) si sa célérité est toujours la même (respectivement varie dans le temps). On peut appliquer ceci à la variation elle-même. Si cette variation est toujours la même, le mouvement est dit uniformément difforme ; il est difformément difforme dans le cas contraire. Et l'on peut itérer ainsi autant que l'on veut. Il n'est pas difficile de voir ici, en germe, la notion d'accélération, puis d'assimiler le mouvement 'uniformément difforme' à notre mouvement 'uniformément accéléré', et ainsi de suite.

Nicole Oresme, un fameux 'Parisien' du 14$^{ème}$ siècle, utilisait même des techniques rappelant pour nous celles du calcul intégral, voire des coordonnées cartésiennes, afin de calculer la 'célérité globale'. Par exemple, celui du mouvement uniformément difforme était donné par le diagramme suivant :

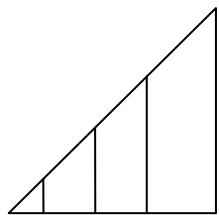

Les segments le long de la droite horizontale, représentant spatial du temps, étaient nommés '*longitudines*', ceux le long de la droite verticale '*latitudines*'.

Cette forme de représentation remonte à la 'méthode des indivisibles', dont, si l'on en croit Archimède, Démocrite se serait déjà servi. Le Syracusain l'utilisait à sont tour, pour calculer certaines surfaces limitées par une parabole. Il paraît donc qu'elle était bien connue des mathématiciens de l'Antiquité. Elle semble toutefois avoir été considérée avec une certaine suspicion, car elle supposait l'existence permanente d'un reste.

L'attitude des 'Oxfordiens' et des 'Parisiens' semble plus détendue, dans la mesure où les seules mathématiques étaient concernées, et ils utilisaient ces calculs beaucoup plus fréquemment. On pourrait en déduire que les principaux obstacles (quantification et utilisation d'une infinité d'instants) pour définir le concept de vitesse (au sens de la physique moderne) avaient été plus ou moins surmontés, par les mathématiciens, bien avant Galilée.

Malheureusement, ceci n'est qu'une reconstruction intellectuelle. Il suffit de rappeler les difficultés à préciser le sens des termes liés au mouvement, ainsi ceux d'*impetus* ou de degrés de vitesse (cf. III.1.γ). Cela reste vrai après Galilée ; ainsi chez Descartes ou même Leibniz, la distinction n'est pas claire, entre ce qu'on appellerait, de nos jours, travail, énergie ou force, ni même entre vitesse et accélération. En employant sans précaution des termes qui n'ont acquis leur précision qu'après un travail lent et considérable, on efface un pan entier de la physique.

Car si les mathématiciens pensaient avoir le droit d'utiliser des artifices (tels les indivisibles), c'est précisément parce qu'ils étudiaient des objets mathématiques (et non physiques). Cela était déjà vrai de l'Antiquité grecque et explique sans doute l'extrême difficulté du statut de l'infini dans le monde fini d'Aristote (*Physique*, III, 5-8 ; VI, 2 ; VI, 9 ; VIII, 8 ; *Du Ciel*, 299a11-17 ; *Métaphysique*, 1071b19 ; cf. §I.3'). Cela le restait à la Renaissance comme le montre l'emploi des 'imaginaires'. De plus, cette situation présentait un avantage pour leur tranquillité. Ainsi le cardinal Bellarmin la garantissait à tous les coperniciens qui s'en tenaient au cadre mathématique (lettre à Paolo Foscarini, XII, p. 171), et on sait les tracas que Galilée dût subir pour s'y être refusé.

Le concept mathématique de 'célérité globale' n'était pas celle d'un mobile physique, moins encore la distance parcourue. Par exemple, dans un cadre mathématique, peu importe si l'on considère la vélocité comme dépendant du temps, de la distance ou des phases de la lune (cf. en particulier, III.1.δ). Par contre, ce sera l'un des problèmes cruciaux à résoudre pour fonder une science du mouvement. Le meilleur argument *a contrario*, est sans doute, le long délai de trois siècles entre les œuvres d'un Oresme par exemple, et les concepts relatifs au calcul infinitésimal. Ainsi que le note Koyré, il s'agissait de rien moins que de 'réformer la structure de notre intelligence elle-même' ([KOY1], p. 171).

## III. Galilée et le concept de vitesse.

### 1. Continuité et mouvement.

La physique galiléenne se fonde sur deux concepts fondamentaux.

**Le premier affirme que tout ce qui a lieu dans la nature est une approximation d'une certaine réalité (ou si l'on préfère d'une théorie) mathématique, puisque la nature est écrite mathématiquement.**
**Le second, conséquence du premier, est la nature continue de la réalité.**

La mathématisation de la nature, fondement de la science galiléenne, est affirmée à plusieurs reprises, et on la trouve déjà dans le *Messager des étoiles*. Le principe de continuité, parce qu'**utilisé** seulement dans les démonstrations, n'est pas formulé explicitement.

Le premier principe est plus large, et en un sens est philosophique, Galilée prenant appui sur Platon. Le second permet de rompre avec les conceptions *péripatéticiennes* (i.e. de l'aristotélisme scolastique). Cette rupture s'exprimant en effet plutôt au travers de la continuité de (ou dans) la nature, que de sa quantification (cf. §I.3).

Suivant Aristote, il y a bien, en général, proportionnalité entre la force motrice et son effet (le mouvement). Par contre, cette proportionnalité échoue si l'on considère des quantités assez petites. Ainsi une force peut mouvoir un corps, alors qu'une plus petite en sera incapable. Un navire ne saurait être mû que par un nombre minimal de haleurs, et un seul homme n'arrivera pas même à l'ébranler ([ARI], VII, 249b20-250a24). En termes modernes, on reconnaît l'exemple d'une discontinuité d'une fonction caractéristique.

Cet aspect fluctuant de la nature aristotélicienne n'attend pas le physicien pisan pour être dénoncé. Ainsi lorsque le platonisant Plutarque rend compte des formidables engins inventés par Archimède, traditionnellement considéré comme un platonicien, il prend précisément l'exemple d'un énorme navire, que le savant syracusain parvient, seul, à mouvoir à son gré ([PLU], Vie de Marcellus, 22, p. 680). Il réfute ainsi à la fois une thèse universellement admise et une thèse d'Aristote, celle de la discontinuité de la cause et de l'effet. En se réclamant d'Archimède, Galilée non seulement affirme la nature mathématique du monde, mais surtout sa continuité.

A. Koyré crédite, à tort nous semble-t-il[3], Descartes et non Galilée, d'être le véritable auteur du principe d'inertie (à savoir 'si aucune force extérieure n'est appliquée à un corps, ce dernier se trouve soit au repos, soit animé d'un mouvement rectiligne uniforme'). Quoiqu'il en soit, Descartes le tient pour exact. Par contre, lorsque le savant français donne ses lois sur les chocs de corps en mouvement, aucune n'est conforme avec le concept de continuité (cf. Annexe). Ainsi que l'affirmera Huygens, en établissant les lois des chocs (élastiques) de la mécanique classique, elles sont pratiquement toutes fausses. Et pour Leibniz, elles le sont précisément en ce qu'elles violent ce principe (cf. [LEI1], VI, p. 131, 249). Très paradoxalement, Descartes, ce pourfendeur sans pitié de la philosophie aristotélicienne, reste certainement plus aristotélicien que Galilée.

Dans le *Dialogue sur les deux grands systèmes du monde*, le principe d'inertie est donné sous la forme de la permanence du mouvement uniforme d'un corps se déplaçant sur un plan horizontal. La démonstration repose sur une suite de lemmes de continuité, ce qu'on appelle le théorème des valeurs intermédiaires : une fonction numérique continue qui prend deux valeurs distinctes (soit *a* et *b*) en deux points distincts (soit *x* et *y*), prend toutes les valeurs intermédiaires (de *a* à *b*) entre ces deux points (*x* et *y*). En termes actuels on dirait que l'image d'un ensemble connexe est connexe. Le principe d'inertie apparaît d'ailleurs comme résultat intermédiaire d'un raisonnement sur la continuité des degrés de vitesse, du repos au mouvement.

La démonstration est faite en plusieurs étapes que nous résumons.

α) Définition de l'égalité de la rapidité.

Salviati rejette la notion usuelle (et aristotélicienne) où l'on compare la rapidité de deux mouvements en fixant soit la distance soit le temps, au profit d'une définition 'plus universelle', dépendant à la fois des distances et des temps. Ainsi 'des vitesses sont dites égales quand les espaces parcourus sont entre eux dans le même rapport que les temps employés à les parcourir' (VII, p. 48). La portée de cette nouvelle **définition** (déjà énoncée par Archimède dans son *Traité des spirales*, mais sous la forme d'une proposition) est considérable, mais par manque de place, nous renvoyons le lecteur intéressé à [OFM3]. C'est elle en effet qui permet de traiter la rapidité comme une grandeur. On ne pourra certes pas comparer n'importe quels mouvements, certains étant à des moments plus rapides, à d'autres moments plus lents. Mais cela signifie que si deux mouvements sont comparables, alors on pourra les comparer en ne considérant que des temps aussi petits que l'on veut. La rapidité devient ainsi une propriété *localisable*.

---
[3] D'autant qu'il réduit, citant Pascal, les énoncés de Galilée sur l'inertie à 'un mot [écrit] à l'aventure sans y faire une réflexion plus longue et plus étendue' ([KOY3], p. 162).

Cette définition devient une proposition dans les *Discours*, s'appuyant sur quatre axiomes (VIII, p. 192-193) et sur l'énoncé V.5 des *Éléments* d'Euclide (cf. §I.3'). L'importance de cette définition V.5 est telle, aux yeux de Galilée, qu'il lui consacrera quelques temps avant sa mort un travail nommé ultérieurement 'la cinquième journée'[4].

β) Le paradoxe de la continuité des degrés de vitesse.

Il s'exprime dans l'étonnement de Sagredo qui 'trouve un peu fort que ce boulet de canon (…) dont la chute, à ce que nous voyons est si précipitée qu'en moins de dix battements de pouls, il franchira plus de cent brasses (57mètres environ), ait passé, au cours de son mouvement, par un degré de vitesse tellement petit que, s'il avait continué à se mouvoir à ce degré de vitesse-là, sans s'accélérer, il n'aurait pas parcouru la même distance en un jour entier'. Et, rajoute aussitôt Salviati, pas même 'en un an, en dix ans en mille ans' (VII, p. 46).

γ) La signification des degrés de vitesses et de lenteur.

Comme l'*impetus* dont elle est quelquefois synonyme chez Galilée, et que l'on peut au mieux décrire comme étant 'une espèce de' (cf. [KOY3], p. 49), le degré de vitesse est une notion assez floue qui essaie de quantifier la célérité, par comparaison à un autre mouvement. On ne considérera donc que des valeurs entières, d'où une arithmétique du mouvement. Quant aux degrés de lenteur, qui jouent un rôle important dans la démonstration, ce sont des compagnons des précédents, en quelque sorte les inverses des degrés de vitesse. La difficulté dans le *Discours*, est l'utilisation par Galilée de ces notions traditionnellement arithmétiques pour exposer sa conception continue de vitesse. D'où une certaine ambiguïté, entre cette nouvelle interprétation défendue par Salviati, et celle plus familière de ses interlocuteurs. Elle s'exprime par exemple dans l'incrédulité de Sagredo (cf. β)).

δ) La propriété fondamentale.

i) Acceptée sans la moindre difficulté, elle énonce que la vitesse acquise par un corps est indépendante de sa trajectoire et ne dépend que de la hauteur de chute. Ceci est vrai aussi bien pour les plans inclinés,

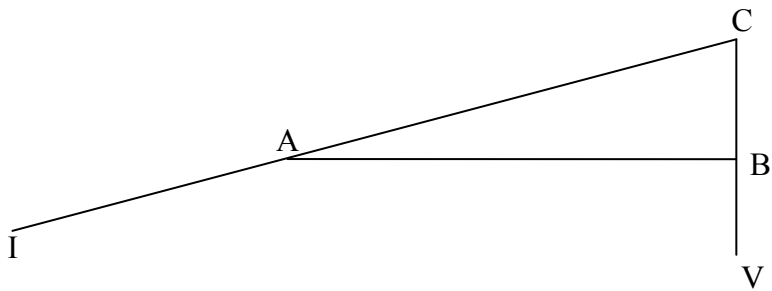

Figure 1

que pour les pendules.

---

[4] C'est Torricelli qui l'a mis en forme à partir d'une dictée de Galilée aveugle, d'où d'importants problèmes pour savoir dans quelle mesure le texte représente bien la pensée de Galilée. On en trouvera une version (en italien) avec diverses variantes dans [GIU] et une traduction anglaise dans [DRA]. En outre Enrico Giusti, dans l'ouvrage précédemment cité, donne une étude détaillée de la reconstruction du livre V par l'école galiléenne.

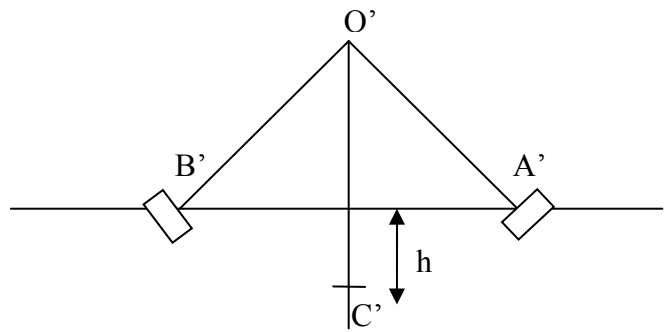

Figure 2

Cette propriété, selon Sagredo, est conséquence évidente de ce que la célérité est fonction de la distance au centre de la terre dont se rapproche le corps pendant sa chute.

ε) La double continuité.

i) Un long raisonnement utilise alors la continuité du mouvement le long de plans inclinés. D'une part, il dépend continûment de l'angle d'inclinaison, de l'horizontale à la verticale, et d'autre part, de la distance au point de départ.

En appliquant la propriété fondamentale δ), et en considérant des plans aussi proches que l'on veut de l'horizontale, Galilée en déduit que le mobile doit passer par tous les degrés de lenteur aussi grands que l'on veut, le degré de lenteur infini correspondant au repos. Il doit donc inversement passer par tous les degrés de vitesse aussi petits que l'on veut, à partir de zéro, degré de vitesse au repos (raisonnement qui sera repris par Leibniz pour poser sa 'Loi de continuité' dans les *Nouveaux essais* ([LEI2], préface)). L'importance de ce point apparaîtra ci-dessous. Si Sagredo agrée, il n'en va pas de même de Simplicio, qui demande une démonstration qualitative.

ii) Toujours sous l'hypothèse 'suivant le cours ordinaire de la nature et tout obstacle extérieur étant écarté', la propriété fondamentale implique le principe d'inertie au sens où, 'sur le plan horizontal, [un mobile] n'acquerra jamais naturellement aucun degré de vitesse quel qu'il soit, attendu que, sur un tel plan, il ne se mettra jamais en mouvement de lui-même', mais 'une fois acquis, il se continuera toujours à une vitesse uniforme' (VII, p. 52-53). La démonstration est la suivante. D'après i), le mouvement d'un mobile le long d'un plan incliné quelconque (mais non horizontal) prend tous les degrés de vitesse. On peut donc considérer que, pour tout corps en mouvement, celui-ci a été acquis, à partir du repos, par descente le long d'un certain plan K'D, éventuellement vertical (cf. fig. 2' ci-dessous).

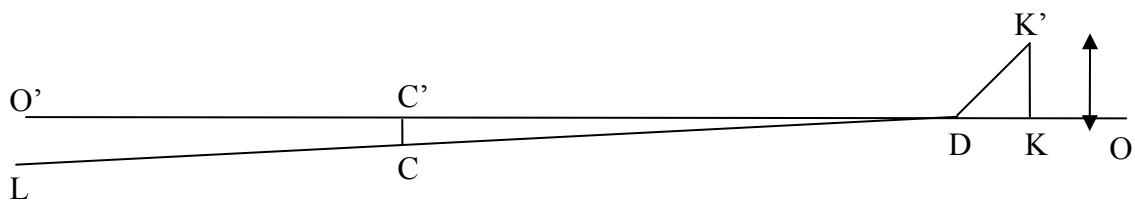

Figure 2'

Soit C' un point arbitraire de l'horizontale OO'. On peut choisir le plan DL d'inclinaison si faible, en sorte que le segment C'C soit très petit, où C désigne la projection verticale de C' sur KL. D'après la propriété fondamentale, la différence de rapidité en C par rapport à C' est mesuré par la distance C'C. En effet la hauteur de chute en C est égale à K'K+C'C. Et puisque C'C est aussi petit que l'on veut la célérité en C est très proche de celle en C', qui elle-même est supérieure à celle en D. Ceci étant vrai pour tout point C' de OO' (le choix du plan DL dépendant toutefois du point C' sur OO'), le mouvement se continue indéfiniment sur l'horizontale, et puisque C'C est très petit, ce mouvement est uniforme.

La démonstration met en évidence les difficultés à surmonter pour arriver à une notion physique permettant de quantifier la célérité en tant qu'objet physique. Ainsi Sagredo confond-il vitesse moyenne et rapidité instantanée ; c'est pour le détromper que Salviati fait cette digression sur le mouvement le long des plans inclinés. La confusion de Sagredo repose sur une conception arithmétique des degrés de vitesse : pour un mouvement partant du repos, il doit exister un plus petit degré de vitesse, surtout si le mouvement est très rapide. Répétons-le, on a là un concept extrêmement flou, élaboré pour rendre compte numériquement d'une qualité (mouvement plus ou moins rapide ou lent). C'est précisément la continuité de la célérité qui va permettre cette mathématisation en éliminant les embarras qui s'étaient accumulés sur cette question (ce 'bourbier' (*quagmire*) dont parle Russell)[5].

Si l'on décide d'identifier degré de vitesse et vitesse (instantanée), le raisonnement apparaît archaïque, difficilement compréhensible et surtout inutile. Au contraire, l'importance du texte tient au passage d'une conception arithmétique de la rapidité (conçue qualitativement et mesurée par des degrés) à une conception continue (conçue comme une grandeur). Ce n'est qu'après cette leçon, qu'une telle identification devient envisageable. Enfin, le texte témoigne d'un traitement quasi infinitésimal du problème, qui suit de la continuité générale des phénomènes, y compris lorsqu'on s'approche de zéro (i.e. du repos dans le cas du mouvement).

## 2. La notion de Vitesse ou qu'est-ce que la vitesse.

Dans cette partie du *Dialogue*, la continuité permet de contourner l'absence de la vitesse (instantanée). La complexité de la démonstration provient des nécessaires allers-retours entre les divers sens que recouvrent la célérité, et que tentent de quantifier les degrés de vitesse. Ceux-ci ont en effet un sens, en tout point du mouvement décrit par un corps, tandis que la vitesse moyenne n'en a que pour des parties finies du mouvement. De la continuité découlent deux conséquences fondamentales qui, pour Galilée, ne sont que les deux aspects d'une même idée.

Du point de vue physique, elle permet de confronter théorie et expérience, celle-ci pouvant être considérée comme conséquence de celle-là. Même si les données ne correspondent pas exactement, il suffit qu'elles en soient suffisamment proches. C'est la très grande innovation que Galilée ne cesse d'utiliser, et qui signe peut-être le mieux la naissance d'une physique mathématique : l'approximation. Et en effet, la nature s'écrirait-elle en termes mathématiques, elle n'en demeurerait pas moins indéchiffrable dans un univers qui ne serait pas mathématique. On aurait alors en effet deux mondes disjoints, suivant le même raisonnement que celui de Parménide, à propos du monde des Idées et du monde sensible ([PLA2], 133b-134c, 134c-135a).

S'il est possible d'étudier la chute des corps, c'est que, quoique le milieu ambiant soit aérien, le désaccord entre le fait et ce que prévoit la théorie, dépend continûment de la résistance qu'oppose le milieu aux corps, résistance qui est nulle dans le vide. Le résultat sera vrai, dès lors qu'on se place dans des conditions suffisamment proches de celles du vide. C'est-à-dire en ne considérant que des corps pour lesquels la résistance de l'air est suffisamment faible (VIII, p. 118-119), et en acceptant une marge d'erreur également faible.

---

[5] Une autre difficulté sur laquelle insiste Galilée, est ce qu'on appellerait aujourd'hui, le caractère continu mais non uniformément continu de la vitesse (cf. VII, p. 47-51).

La preuve ne se trouve pas dans l'exactitude du résultat prévu par la théorie, mais dans le rapport prédictif entre celles de Galilée et d'Aristote (cf. VIII, p. 120).

Ainsi Galilée peut-il affirmer simultanément que l'arrivée de deux corps de poids différents lâchés d'une certaine hauteur n'est pas concomitante, et qu'ils arrivent au même moment.

'Je ne voudrais pas, signor Simplicio, qu'à l'exemple de tant d'autres, détournant notre propos de notre objet principal, vous vous attachiez à telle chose que j'ai dite et qui s'écartait de la vérité de l'épaisseur d'un cheveu, pour essayer de dissimuler, sous ce cheveu, l'erreur, aussi grosse qu'une amarre, qu'un autre a commise. Aristote dit : « Une boule de fer de cent livres, tombant d'une hauteur de cent brasses, arrive au sol avant qu'une boule d'une livre soit descendue d'une seule brasse ». Je dis, moi, qu'elles arrivent en même temps.' (VIII, p. 109).

C'est la continuité qui prouve à la fois la justesse de sa théorie et l'erreur d'Aristote. Sans elle, les deux théories seraient fausses, l'une l'étant seulement (beaucoup) moins que l'autre.

D'un point de vue méta-physique (au sens de la méta-mathématique, mais également au sens philosophique), elle permet une double remise en cause de la pensée péripatéticienne. D'une part, en montrant l'impossibilité de se fier imprudemment à l'évidence. D'autre part, en mettant en évidence les erreurs d'une théorie dont le seul objectif est de rendre compte de l'expérience sensible, mais aussi, et pour Galilée c'est le plus important, en expliquant les causes de ces erreurs.

Ce n'est qu'après de nombreuses années d'étude et d'expériences, que Galilée établissait une loi (mathématique, physique ?) mettant en relation la distance couverte par un mobile en chute libre et le temps de parcours. Auparavant, il avait conclu, contre les thèses aristotéliciennes, à la composition des divers mouvements (naturels et violents) et surtout à l'éternité du mouvement (rectiligne) uniforme (cf. paragraphe précédent).

La simplicité des moyens utilisés dans ces expériences peut surprendre : plans inclinés et pendules (cf. VII, p. 44-55 ; VIII, p. 126-129). Galilée souligne l'importance des pendules et aussi l'originalité (au moins par rapport aux 'philosophes péripatéticiens') de leur étude (cf. VIII, p. 139). Et pourtant, il fait dire à Sagredo que cette simplicité même pourrait être source de railleries contre lui, Galilée (VIII, p. 131, p. 140). Mais en outre, les concepts mathématiques utilisés ne sont guère plus spécialisés, et ne vont certainement pas au-delà des *Eléments* d'Euclide.

Le phénomène des pendules une fois compris (comme **une** forme de chute d'un corps), son utilisation devenait une nécessité absolue pour minimiser la friction :

'Allant plus loin, j'ai voulu éliminer les inconvénients qui pouvaient naître du contact des mobiles avec le plan incliné ; et finalement (…) j'ai attaché [deux balles, l'une de plomb, l'autre de liège à deux cordons très fins et] écarté les deux balles de la perpendiculaire' (VIII, p. 128).

La symétrie du mouvement du pendule, à la fois exprimait et était interprétée par l'éternité du mouvement, quoique d'un mouvement non uniforme (VIII, p. 129, 138 ; et aussi paragraphe précédent).

Mais, résultat complémentaire, Galilée obtenait une nouvelle loi, sorte 'd'équivalence' entre la hauteur de chute d'un corps et sa vitesse. La vitesse d'un corps tombant d'une certaine hauteur est indépendante de ce corps (sa forme, son poids…) et de la trajectoire (chute libre, plan incliné ou trajectoire circulaire d'un pendule) (VIII, p. 128 et 1.δ).

Cette loi, que nous avons nommée la propriété fondamentale (de la chute des corps) n'est peut-être pas aussi célébrée que d'autres travaux de Galilée, particulièrement ceux concernant l'astronomie et le système copernicien. Ainsi Geymonat voit-il dans l'invention de 'sa fameuse lunette (…) l'origine de la phase la plus glorieuse de son activité scientifique' ([GEY], p. 44). Déjà du vivant de Galilée, les éditeurs des *Discours*, Louis et Abraham Elzévir, pour vanter les mérites de Galilée, mettaient en avant ses découvertes faites 'grâce à la lunette astronomique' (VIII, p. 48). Ce résultat est pourtant crucial pour élaborer le concept de vitesse et, en conséquence, le calcul infinitésimal. Déjà L. Brunschvicg soulignait, il est vrai suivant un point de vue différent, l'importance des travaux de Galilée dans la 'découverte' du calcul infinitésimal, mais surtout intégral ([BRU], 10, p. 215-6).

Naturellement, il est possible de donner, à ce que nous avons nommé 'célérité', le sens de vitesse moyenne déjà présente chez Aristote (cf. I.3'). Mais comme on va le voir, se servir de cette interprétation, pour sauver ce qui peut l'être de la physique aristotélicienne du mouvement, est extrêmement trompeur.

Si l'on se réfère à la 'célérité', le problème fondamental était, suivant notre analyse, le manque d'unité de mesure (cf. §I.3' de cette partie). Mais c'était aussi, dans un autre sens, le manque d'unité par impossibilité de lui donner un sens unique global[6].

Sur ce point, les scolastiques n'ont guère apporté de changement. Malgré les résultats obtenus par l'introduction des différents mouvements uniforme, uniformément difforme, uniformément difformément difforme…, il restait à caractériser un mouvement comme appartenant à l'une ou l'autre de ces classes. Ainsi, dans la classification de Nicole Oresme, qui pour certains préfigurent les travaux de Galilée (cf. §II.2), les différents cas étaient disjoints les uns des autres, chacun demandant une étude spécifique, d'où l'impossibilité à définir **la** vitesse.

Car une telle définition quantitative (ou mathématique, ou physique), pour être indépendante du mouvement, doit être prise au sens **de vitesse instantanée**. C'est pourquoi, la question de la vitesse n'était jamais physique, mais mathématique. Aussi avons-nous préféré utiliser le terme de 'célérité' et non de 'vitesse' tant qu'une telle définition était proprement impossible, alors même que la distinction qualité/quantité s'estompait au Moyen-âge, si jamais elle eût vraiment une grande importance dans la Grèce antique (cf. chap. I).

Il reste à comprendre pourquoi, avec les travaux de Galilée, cette notion prenait sens, mais également en quoi la propriété fondamentale de la chute des corps (cf. 1.δ) modifiait radicalement la conception antérieure du mouvement.

La réponse est étonnamment simple : cette propriété fondamentale permettait de **calculer** la vitesse (instantanée).

En effet, la continuité des degrés de vitesse, conçus par Galilée, contre la tradition, comme des grandeurs, implique que tout mouvement instantané se ramène à une chute verticale suivie d'un mouvement uniforme. Il en est de même du mouvement des planètes, et le savant italien se trouve si assuré de l'exactitude de cette théorie, qu'il affirme pouvoir inverser le calcul : de la connaissance précise de leur mouvement de rotation, on peut obtenir la hauteur de leur chute.

'Imaginons que le divin Architecte (…) ait placé, immobile en son centre, le Soleil (…) Demandons-nous maintenant à quelle hauteur à quelle distance du soleil était le lieu où ces globes avaient été primitivement créés. Pour résoudre ce problème, il faut nous enquérir, auprès des plus habiles astronomes, des grandeurs des cercles suivant lesquels tournent les planètes ainsi que leur temps de leurs révolutions.' (VII, p. 53).

---

[6] Qui n'est rien d'autre que la question de Socrate dans le *Théétète*, ou selon sa reformulation par Théétète, la difficulté à 'ranger [une multitude] dans une unité' (147d).

On peut considérer cela comme 'une fable platonicienne' de la création du monde (cf. par exemple l'introduction à la traduction française par R. Fréreux du *Dialogue*). Et pourtant, c'est ainsi que Newton va raisonner pour obtenir le mouvement des planètes autour du soleil[7].

**Dès lors, et ce quel que soit le mouvement et la trajectoire, uniforme ou difforme, rectiligne ou circulaire, il suffit, pour obtenir le degré de vitesse *v* d'un mobile au temps *t*, de calculer la hauteur *H* que ce mobile est capable d'atteindre. Ou bien, inversement, connaissant ses 'degrés de vitesse', il est possible d'en déduire la hauteur *H*, ainsi dans le cas des planètes.**

Les conséquences suivantes découlent immédiatement des équations que donnent Galilée pour un mouvement de chute libre. Le degré de vitesse du mobile au temps *t* lui permettrait, cette vitesse restant constante, de franchir une distance égale à 2 fois celle qu'il a parcourue depuis le repos.
On obtient ce résultat par une sorte d'intégration, à la façon de l'exemple du §II.2 (cf. VII, p. 52 et VIII, p. 243, scolie du problème 9).
On aura donc : $2H = vt$ (*).
D'après la loi de la distance dépendant du carré du temps (les constantes pouvant toujours être supposées égales à 1), on a : $H = t^2$,
et de l'égalité (*) on déduit :
$v^2 = 4H$ (1)
ou encore
$H = v^2/4$ (1')

Ces résultats, valables *a priori* pour la seule chute libre, sont (d'après 1.ϵ.i) universellement vrais. En effet, dès lors que l'on considère les degrés de vitesse comme des grandeurs continues, quelque soit le mouvement qu'elles permettent de mesurer (i.e. d'en mesurer la rapidité) à un temps donné, elles sont égales à celles d'un mouvement en chute libre à un temps donné. En d'autres termes, la mesure de la rapidité d'un mouvement est **indépendante de la nature de ce mouvement**. C'est l'équation (1).
D'après cette équation, la grandeur *v* est déterminée par la racine carrée de la hauteur H. **Autrement dit, la hauteur donne un étalon (une unité) universel pour calculer les 'degrés de vitesse' galiléens, c'est-à-dire la vitesse (car définie pour tout temps, tout mouvement et toute trajectoire).**

Aussi surprenant que cela paraisse, on pouvait obtenir et mesurer la vitesse (instantanée) avant qu'elle ait été définie.

Si l'égalité (1) unifie le mouvement, en rendant possible une mesure universelle de la vitesse (indépendante du mouvement), l'égalité (1'), en unifiant les 'forces' (au sens flou qu'on trouve encore chez Galilée) dues à la vitesse et à la position, permet d'obtenir une définition universelle de ce qu'on appelle maintenant 'énergie'. Elle exprime essentiellement la **conservation de l'énergie**[8], i.e. de la somme énergie potentielle et énergie cinétique[9]. On pourrait objecter que les masses n'apparaissent pas. Mais parce que l'on étudie ici le mouvement d'un seul corps de masse constante, on peut la supposer égale à l'unité.

---
[7] À la différence toutefois de Galilée, il ne supposera pas que le mouvement est nécessairement circulaire, mais montrera, précisément par ce raisonnement, qu'il doit être elliptique.
[8] 'la loi galiléenne de la chute des corps conduisait infailliblement à ce résultat [en termes modernes, la conservation de l'énergie cinétique]' ([GUE], p. 62), et même un peu plus.

Ces équations sont au fondement de ce qu'on appellera la querelle des 'forces vives'. Elle portait sur ce qui se conservait dans le mouvement, mettant aux prises cartésiens et leibniziens. Pour les premiers, c'était la vitesse (en tant que grandeur), pour les seconds, son carré, les leibniziens se trouvant alors du côté de Galilée. Cette fois encore, Descartes est plus proche des péripatéticiens, ou plus éloigné de la physique moderne, que Galilée[10].

### 3. Vitesse et infinitésimaux.

Ainsi, c'est de la physique de Galilée (et seulement à partir d'elle), qu'il devenait possible de calculer la vitesse (instantanée). Et cela signifiait simultanément **la possibilité de la définir**. Si les philosophes/physiciens précédemment n'avaient pas été capables de mesurer la vitesse, ce n'était pas dû à des blocages psychologiques, induits par la théorie des catégories aristotéliciennes au travers de l'opposition qualité et quantité. Au contraire, la théorie d'Aristote était conçue (au moins en partie) pour répondre aux difficiles problèmes que posait le mouvement, dont les paradoxes de Zénon. L'impossibilité de mesurer la 'célérité' d'un mobile ne tenait pas à quelque erreur concernant son appartenance à telle catégorie, mais parce qu'il n'y avait rien à mesurer, seulement une infinité de vitesses différentes associées à une infinité de mouvements différents causés par une infinité de moteurs distincts. Et Aristote avait raison de considérer alors seuls les temps et distances comme mesurables.

Mais qu'en est-il des critiques portant sur son incompréhension de la 'vitesse moyenne', laquelle est définissable dès que l'on conçoit une mesure de la distance et du temps ? C'est négliger qu'elle n'a de sens qu'à la condition d'ajouter l'intervalle de temps pendant lequel le mouvement a lieu, à cette vitesse (moyenne). Sinon en effet, la tortue serait plus rapide qu'Achille, pourvu qu'elle se meuve suffisamment longtemps. On a d'ailleurs le même paradoxe pour le mouvement uniformément accéléré qui n'est pas uniformément continu (cf. III.1.δ et n. 5 §III.1 ; VII, p. 47-51), et on le retrouve encore à la fin du 18$^{ème}$ siècle dans la querelle sur les 'forces vives' ([GUE], p. 142).

Mais, connaître le temps **et** la 'vitesse moyenne', est exactement équivalent à connaître le temps et la distance. Et c'est bien ce dernier couple qui est réellement mesuré, aucun des deux d'ailleurs ne donnant d'indication sur la vitesse (i.e. instantanée) du mobile, comme le montre le paradoxe d'Achille et de la tortue. En outre, la notion de vitesse moyenne, en tant que grandeur, n'était certainement pas étrangère au Stagirite, et donc *a fortiori*, aux mathématiciens de la Grèce ancienne (cf. §I.3').

Finalement, il ne s'agit plus de s'interroger sur les obstacles à une mathématisation du mouvement, mais de comprendre comment l'étude d'un phénomène extrêmement restreint, celui de la chute des corps, a conduit à ce formidable bouleversement que furent les infinitésimaux. Ils vont surgir en effet de travaux sur un mouvement très particulier, celui induit par la gravité, c'est-à-dire en termes mathématiques, d'un mouvement rectiligne dont la vitesse varie uniformément.

La nouveauté fondamentale était la suivante : dès lors que le concept de vitesse prenait sens, i.e. qu'elle pouvait être mesurée (en tant qu'instantanée) par son équivalence à un espace, la vitesse ne dépendait plus du mouvement ou du moteur, et devenait universelle. Et Galilée, presque malgré lui, devenait 'un des inspirateurs de cette géométrie supérieure d'où

---

[9] C'est l'argument de Leibniz pour rejeter la loi cartésienne de conservation du mouvement (cité in [GUE], p. 65, n. 2).
[10] Pour un point de vue contraire, cf. [GUE], p. 63 par exemple. La position d'Alexandre Koyré est plus ambiguë.

devait sortir l'application du calcul infinitésimal aux réalités dynamiques' ([GUE], p. 82, n. 2).

De nouveaux problèmes apparaissaient, les problèmes inverses. Il s'agissait de trouver un mouvement ou une trajectoire dont la vitesse avait certaines propriétés (cf. paragraphe précédent, l'exemple cosmogonique donné par Galilée dans le *Dialogue* (VII, p. 53)). En termes modernes cela consistait à obtenir des solutions d'équations différentielles. Ainsi, Huygens étend-il la loi de constance des pendules (circulaires) pour les petits angles au cas général (cycloïdes), permettant de bâtir des horloges marines qui à leur tour servent au calcul très complexe de la longitude. Et bien entendu en astronomie, cela a conduit à une explication du mouvement des planètes, et finalement à une théorie de la gravitation universelle[11].

## 4. Conclusion

En ce sens, Newton pouvait sans doute affirmer qu'il n'avait fait que mettre ses pas dans ceux de Galilée, un de ces géants sur les épaules desquels il avait trouvé refuge. Néanmoins, ces pas restaient à faire, pour unifier complètement la multiplicité des mouvements.

Si ce problème est bien un problème platonicien, peut-être le problème de Platon par excellence (cf. [OFM2], §III.3 et 4), quelle position adopter dans cette guerre de géants dont nous parle le philosophe athénien ?

Faut-il proclamer Galilée la 'revanche de Platon (…) qui a renversé l'empire d'Aristote' ([KOY2], p. 321), au prix certes d'une alliance diabolique avec Démocrite ?

Faut-il au contraire, selon les mots de M. Clavelin (*op. cité*, p. ix) rapportant ce point de vue pour le critiquer, voir en Galilée 'le dernier disciple d'Aristote' (position adoptée par P. Duhem et W. A. Wallace) ?

Ou encore ne faudrait-il pas s'interroger si, comme dans le *Sophiste* (cf. [OFM2], §III.5), ces deux thèses n'aboutissent pas à leur propre ruine ?

Car, suivant notre analyse, les textes d'Aristote et de Galilée ne diffèrent pas seulement par leur langue, mais par cela même qui est visé sous les termes $\tau\alpha\chi\acute{u}\varsigma$ et *velocitas/velocità*. La première position, loin d'accorder toute son importance aux travaux galiléens sur le mouvement, nie la révolution conceptuelle opérée à leur suite. La seconde, voulant rendre justice à la scolastique héritière d'Aristote, en fait le responsable de ce qu'elle n'eût pas été opérée bien auparavant.

Concevoir la constitution de la vitesse sous une forme épistémologique, revient à s'enfermer dans une telle alternative impossible. Car c'est au contraire sa formalisation mathématique, au travers de la théorie des infinitésimaux, qui lui donne sens. Mais inversement la notion mathématique des infinitésimaux ne se constitue qu'à la suite des travaux des physiciens. L'épistémologie doit faire place ou plus exactement faire une place à l'ontologie. Pour que la vitesse soit connaissable, encore doit-elle être physiquement, c'est-à-dire, être objet de mesure. Du point de vue mathématique, elle est un rapport d'infinitésimaux. C'est à cette unité, être physique/être mathématique qu'aboutit, non sans mal, l'intelligence en lutte avec elle-même.

## *TEXTES CITÉS*

[ARI]   :   Aristote, *Physique,* trad. H. Carteron-L. Robin, Belles Lettres, 1931, 1986

---

[11] 'de ces pages capitales [les discours] est venue l'impulsion qui en moins de cinquante années a permis la mécanisation du système du Monde' dit Maurice Clavelin dans son introduction à la traduction française du *Dialogue* (p. xxii-xxiii).

## *ANNEXE*

## Les lois du choc selon Descartes

'On considère deux corps B et C qui se dirigent l'un vers l'autre, l'un d'entre eux pouvant toutefois être au repos. On note encore B et C leurs grandeurs respectives [c'est une notion assez vague, qu'on peut assimiler de manière très approximative à la masse du corps], V et W leurs vitesses avant le choc, V' et W' après le choc (la vitesse est considérée comme un nombre positif). Bien que cela ne soit pas explicitement précisé par Descartes, leurs trajectoires sont situées sur une même droite. On nomme détermination la direction (et le sens) du mouvement. Ainsi les vitesses de deux mobiles peuvent être égales et avoir une détermination opposée, auquel cas ils se déplacent en sens contraire. Les règles doivent permettre de connaître V' et W' ainsi que leurs déterminations. (…)

Règle 1
Si $B = C$ et $V = W$ alors $V' = V$ et $W' = W$ en changeant tout deux de détermination.
Règle 2
Si $B > C$ et $V = W$ alors $V' = V$ et $W' = W$, et W' change de détermination (seul le texte latin donne toutefois explicitement l'égalité des vitesses après le choc).
Règle 3
Si $B = C$ sont tous deux en mouvement et $V > W$ alors $V' = W' = (V+W)/2$ et W seul change de détermination.
Règle 4
Si $B < C$ et C est au repos, alors $V' = V$ et change de détermination tandis que C reste au repos.
Règle 5
C'est la formule la plus complexe.
Si $B > C$ et C est au repos alors $V' = W' = V-(C/(B+C))V$, V et V' ayant même détermination. Par exemple si B est 2 fois plus grand que C, après le choc B perd le tiers de sa vitesse ; s'il est 3 fois plus grand, il en perd le quart.
( …)
[Règle 6]
La règle 6 porte sur choc entre deux corps B et C, de grandeurs exactement égales, le premier étant en mouvement vers le second, en repos. Pour Descartes, B rejaillirait (son sens s'inverserait) tandis que C serait repoussé en sens opposé.
(…)
Si $B = C$ et C est en repos, alors $V' = (3/4) V$, le sens de V' est opposé à celui de V et $W' = V/4$ est de même sens que V (cf. supra, fig. 1 et 2).' ([OFM1], §30, p. 367-375).